\documentclass[10pt]{amsart}
\usepackage{amsfonts}
\usepackage{amsmath}
\usepackage{amsthm}
\usepackage{amssymb}
\usepackage{latexsym}

\newtheorem{lem}{Lemma}[section]

\newtheorem{prop}{Proposition}[section]

\theoremstyle{definition}
\newtheorem{defn}{Definition}[section]
\theoremstyle{definition}
\newtheorem{thm}{Theorem}

\newenvironment{pf}{\proof}{\endproof}

\theoremstyle{remark}

\numberwithin{equation}{section}
\begin{document}

\newcommand{\thmref}[1]{Theorem~\ref{#1}}
\newcommand{\secref}[1]{Sect.~\ref{#1}}
\newcommand{\lemref}[1]{Lemma~\ref{#1}}
\newcommand{\propref}[1]{Proposition~\ref{#1}}
\newcommand{\corref}[1]{Corollary~\ref{#1}}
\newcommand{\remref}[1]{Remark~\ref{#1}}
\newcommand{\nc}{\newcommand}
\newcommand{\rnc}{\renewcommand}
\nc{\cal}{\mathcal}
\nc{\goth}{\mathfrak}
\rnc{\bold}{\mathbf}
\renewcommand{\frak}{\mathfrak}
\renewcommand{\Bbb}{\mathbb}

\nc{\Cal}{\mathcal}
\nc{\Xp}[1]{X^+(#1)}
\nc{\Xm}[1]{X^-(#1)}
\nc{\on}{\operatorname}
\nc{\ch}{\mbox{ch}}
\nc{\Z}{{\bold Z}}
\nc{\J}{{\cal J}}
\nc{\C}{{\bold C}}
\nc{\Q}{{\bold Q}}
\renewcommand{\P}{{\cal P}}
\nc{\N}{{\Bbb N}}
\nc\beq{\begin{equation}}
\nc\enq{\end{equation}}
\nc\lan{\langle}
\nc\ran{\rangle}
\nc\bsl{\backslash}
\nc\mto{\mapsto}
\nc\lra{\leftrightarrow}
\nc\hra{\hookrightarrow}
\nc\sm{\smallmatrix}
\nc\esm{\endsmallmatrix}
\nc\sub{\subset}
\nc\ti{\tilde}
\nc\nl{\newline}
\nc\fra{\frac}
\nc\und{\underline}
\nc\ov{\overline}
\nc\ot{\otimes}
\nc\bbq{\bar{\bq}_l}
\nc\bcc{\thickfracwithdelims[]\thickness0}
\nc\ad{\text{\rm ad}}
\nc\Ad{\text{\rm Ad}}
\nc\Hom{\text{\rm Hom}}
\nc\End{\text{\rm End}}
\nc\Ind{\text{\rm Ind}}
\nc\Res{\text{\rm Res}}
\nc\Ker{\text{\rm Ker}}
\rnc\Im{\text{Im}}
\nc\sgn{\text{\rm sgn}}
\nc\tr{\text{\rm tr}}
\nc\Tr{\text{\rm Tr}}
\nc\supp{\text{\rm supp}}
\nc\card{\text{\rm card}}
\nc\bst{{}^\bigstar\!}
\nc\he{\heartsuit}
\nc\clu{\clubsuit}
\nc\spa{\spadesuit}
\nc\di{\diamond}

\nc\al{\alpha}
\nc\bet{\beta}
\nc\ga{\gamma}
\nc\de{\delta}
\nc\ep{\epsilon}
\nc\io{\iota}
\nc\om{\omega}
\nc\si{\sigma}
\rnc\th{\theta}
\nc\ka{\kappa}
\nc\la{\lambda}
\nc\ze{\zeta}

\nc\vp{\varpi}
\nc\vt{\vartheta}
\nc\vr{\varrho}

\nc\Ga{\Gamma}
\nc\De{\Delta}
\nc\Om{\Omega}
\nc\Si{\Sigma}
\nc\Th{\Theta}
\nc\La{\Lambda}
\nc\boa{\bold a}
\nc\bob{\bold b}
\nc\boc{\bold c}
\nc\bod{\bold d}
\nc\boe{\bold e}
\nc\bof{\bold f}
\nc\bog{\bold g}
\nc\boh{\bold h}
\nc\boi{\bold i}
\nc\boj{\bold j}
\nc\bok{\bold k}
\nc\bol{\bold l}
\nc\bom{\bold m}
\nc\bon{\bold n}
\nc\boo{\bold o}
\nc\bop{\bold p}
\nc\boq{\bold q}
\nc\bor{\bold r}
\nc\bos{\bold s}
\nc\bou{\bold u}
\nc\bov{\bold v}
\nc\bow{\bold w}
\nc\boz{\bold z}
\nc\bK{\bold K}

\nc\ba{\bold A}
\nc\bb{\bold B}
\nc\bc{\bold C}
\nc\bd{\bold D}
\nc\be{\bold E}
\nc\bg{\bold G}
\nc\bh{\bold h}
\nc\bH{\bold H}

\nc\bi{\bold I}
\nc\bj{\bold J}
\nc\bk{\bold K}
\nc\bl{\bold L}
\nc\bm{\bold M}
\nc\bn{\bold N}
\nc\bo{\bold O}
\nc\bp{\bold P}
\nc\bq{\bold Q}
\nc\br{\bold R}
\nc\bs{\bold S}
\nc\bt{\bold T}
\nc\bu{\bold U}
\nc\bv{\bold v}
\nc\bV{\bold V}

\nc\bw{\bold w}
\nc\bz{\bold Z}
\nc\bx{\bold x}
\nc\bX{\bold X}
\nc\blambda{{\mbox{\boldmath $\Lambda$}}}
\nc\bpi{{\mbox{\boldmath $\pi$}}}
\nc\bpsi{{\mbox{\boldmath $\psi$}}}

\nc\e[1]{E_{#1}}
\nc\ei[1]{E_{\delta - \alpha_{#1}}}
\nc\esi[1]{E_{s \delta - \alpha_{#1}}}
\nc\eri[1]{E_{r \delta - \alpha_{#1}}}
\nc\ed[2][]{E_{#1 \delta,#2}}
\nc\ekd[1]{E_{k \delta,#1}}
\nc\emd[1]{E_{m \delta,#1}}
\nc\erd[1]{E_{r \delta,#1}}

\nc\ef[1]{F_{#1}}
\nc\efi[1]{F_{\delta - \alpha_{#1}}}
\nc\efsi[1]{F_{s \delta - \alpha_{#1}}}
\nc\efri[1]{F_{r \delta - \alpha_{#1}}}
\nc\efd[2][]{F_{#1 \delta,#2}}
\nc\efkd[1]{F_{k \delta,#1}}
\nc\efmd[1]{F_{m \delta,#1}}
\nc\efrd[1]{F_{r \delta,#1}}
\nc{\ug}{\bu^{\text{fin}}}

\nc\fa{\frak a}
\nc\fb{\frak b}
\nc\fc{\frak c}
\nc\fd{\frak d}
\nc\fe{\frak e}
\nc\ff{\frak f}
\nc\fg{\frak g}
\nc\fh{\frak h}
\nc\fj{\frak j}
\nc\fk{\frak k}
\nc\fl{\frak l}
\nc\fm{\frak m}
\nc\fn{\frak n}
\nc\fo{\frak o}
\nc\fp{\frak p}
\nc\fq{\frak q}
\nc\fr{\frak r}
\nc\fs{\frak s}
\nc\ft{\frak t}
\nc\fu{\frak u}
\nc\fv{\frak v}
\nc\fz{\frak z}
\nc\fx{\frak x}
\nc\fy{\frak y}

\nc\fA{\frak A}
\nc\fB{\frak B}
\nc\fC{\frak C}
\nc\fD{\frak D}
\nc\fE{\frak E}
\nc\fF{\frak F}
\nc\fG{\frak G}
\nc\fH{\frak H}
\nc\fJ{\frak J}
\nc\fK{\frak K}
\nc\fL{\frak L}
\nc\fM{\frak M}
\nc\fN{\frak N}
\nc\fO{\frak O}
\nc\fP{\frak P}
\nc\fQ{\frak Q}
\nc\fR{\frak R}
\nc\fS{\frak S}
\nc\fT{\frak T}
\nc\fU{\frak U}
\nc\fV{\frak V}
\nc\fZ{\frak Z}
\nc\fX{\frak X}
\nc\fY{\frak Y}
\nc\tfi{\ti{\Phi}}
\nc\bF{\bold F}

\nc\ua{\bold U_\A}

%%%%%%%%%%%%%%%%%%%%%%%%%%%%%%%%%%%%%%%%%%%%%%%%%%%%%%
\nc\qinti[1]{[#1]_i}
\nc\q[1]{[#1]_q}
\nc\xpm[2]{E_{#2 \delta \pm \alpha_#1}}  %\xpm{j}{l}
\nc\xmp[2]{E_{#2 \delta \mp \alpha_#1}}
\nc\xp[2]{E_{#2 \delta + \alpha_{#1}}}
\nc\xm[2]{E_{#2 \delta - \alpha_{#1}}}
\nc\hik{\ed{k}{i}}
\nc\hjl{\ed{l}{j}}
\nc\qcoeff[3]{\left[ \begin{smallmatrix} {#1}& \\ {#2}& \end{smallmatrix}
\negthickspace \right]_{#3}}
\nc\qi{q}
\nc\qj{q}

\nc\ufdm{{_\ca\bu}_{\rm fd}^{\le 0}}

%%%%%%%%%%%%%%%%%%%%%%%%%%%%%%%%%%%%%%%%%%%%%%%%%%%%%%

%\nc\rtimes
\nc\isom{\cong} 

\nc{\pone}{{\Bbb C}{\Bbb P}^1}
\nc{\pa}{\partial}
\def\H{\cal H}
\def\L{\cal L}
\nc{\F}{{\cal F}}
\nc{\Sym}{{\goth S}}
\nc{\A}{{\cal A}}
\nc{\arr}{\rightarrow}
\nc{\larr}{\longrightarrow}

\nc{\ri}{\rangle}
\nc{\lef}{\langle}
\nc{\W}{{\cal W}}
\nc{\uqatwoatone}{{U_{q,1}}(\su)}
\nc{\uqtwo}{U_q(\goth{sl}_2)}
\nc{\dij}{\delta_{ij}}
\nc{\divei}{E_{\alpha_i}^{(n)}}
\nc{\divfi}{F_{\alpha_i}^{(n)}}
\nc{\Lzero}{\Lambda_0}
\nc{\Lone}{\Lambda_1}
\nc{\ve}{\varepsilon}
\nc{\phioneminusi}{\Phi^{(1-i,i)}}
\nc{\phioneminusistar}{\Phi^{* (1-i,i)}}
\nc{\phii}{\Phi^{(i,1-i)}}
\nc{\Li}{\Lambda_i}
\nc{\Loneminusi}{\Lambda_{1-i}}
\nc{\vtimesz}{v_\ve \otimes z^m}

\nc{\asltwo}{\widehat{\goth{sl}_2}}
\nc\eh{\frak h^e}  
\nc\loopg{L(\frak g)}  
\nc\eloopg{L^e(\frak g)} 
\nc\ebu{\bu^e} 

\nc\teb{\tilde E_\boc}
\nc\tebp{\tilde E_{\boc'}}

\title{Integrable and   Weyl modules \\
for quantum affine $sl_2$. }
%\dedicatory{}
\author{Vyjayanthi Chari}
\address{Vyjayanthi Chari, University of California, Riverside}
\author{Andrew Pressley}
\address{Andrew Pressley, Kings College, London}

\pagestyle{plain} \maketitle
\section{Introduction} Let $\frak t$ be an arbitrary symmetrizable Kac-Moody Lie 
algebra and $\bu_q(\frak t)$ the corresponding quantized enveloping algebra of 
$\frak t$ defined over $\bc(q)$.
If $\mu$ is a dominant integral weight of $\frak t$ then one can associate to it 
in a natural way an irreducible integrable $\bu_q(\frak t)$-module $L(\mu)$. 
These modules have many nice properties  and are well understood, \cite{K}, \cite{L}. 

More generally, given any integral weight $\lambda$,  Kashiwara \cite{K} defined 
an integrable $\bu_q(\frak t)$-module $V^{max}(\lambda)$ generated by an 
extremal vector $v_\lambda$. If $w$ is any element of the Weyl group $W$ of 
$\frak t$, then one has $V^{max}(\lambda)\cong V^{max}(w\lambda)$. Further, if 
$\lambda$ is in the Tits cone, then 
$V^{max}(\lambda)\cong L(w_0\lambda)$, where $w_0\in W$ is such that 
$w_0\lambda$ is dominant integral. In the case when $\lambda$ is not in the Tits cone, the module $V^{max}(\lambda)$ is not irreducible and very 
little is known about it, although it is known that it admits a crystal basis, 
\cite{K}.

In  the case when $\frak t$ is an affine Lie algebra, an integral weight 
$\lambda$ is not in the Tits cone if and only if $\lambda$ has level zero.  Choose 
$w_0\in W$ so that $w_0\lambda$ is dominant with respect to the underlying 
finite-dimensional simple Lie algebra of $\frak t$.  In as yet unpublished work, 
Kashiwara proves that $V^{max}(\lambda)\cong W_q(w_0\lambda)$, where 
$W_q(w_0\lambda)$ is an integrable $\bu_q(\frak t)$-module  defined by 
generators and relations analogous to the definition of $L(\mu)$. 

In \cite{CP5}, we studied the modules $W_q(\lambda)$ further. In particular, we 
showed that they have a family  $W_q(\bpi)$  of non--isomorphic finite-dimensional quotients  
which are maximal, in the sense that any another finite-dimensional quotient is 
a proper  quotient of some $W_q(\bpi)$. In this paper, we show that, if $\frak t$ is the 
affine Lie algebra associated to $sl_2$ and $\lambda=m\in\bz^+$, the modules 
$W_q(\bpi)$ all have the same dimension $2^m$. This is done by showing that the 
modules $W_q(\bpi)$, under suitable conditions, have a $q=1$ limit, which allows 
us to reduce to the study of the corresponding problem in the classical case 
carried out in \cite{CP5}. The  modules $W_q(\bpi)$ have a unique irreducible 
quotient $V_q(\bpi)$, and we show that these are all the irreducible 
finite-dimensional $\bu_q(\frak t)$-modules.  In \cite{CPqa}, \cite{banff},  a 
similar classification was obtained by regarding $q$ as a complex number and 
$\bu_q(\frak t)$ as an algebra over $\bc$; in the present situation, we have to 
allow modules defined over finite extensions of $\bc(q)$.   

We are then able to realize the modules $W_q(m)$ as being the space of 
invariants of the action of the Hecke algebra $\cal{H}_m$ on the tensor product 
$(V\otimes\bc(q)[t,t^{-1}])^{\otimes m}$, where $V$ is a two-dimensional vector 
space over $\bc(q)$. Again, this is done by reducing to the case of $q=1$. 

In the last section, we indicate how to extend some of the results of this paper 
to the general case. We conjecture that the dimension of the modules $W_q(\bpi)$ 
depends only on $\lambda$, and we give a formula for this dimension.

\section {Preliminaries and Some Identities} 
Let  $sl_2$ be the complex Lie algebra with basis $\{x^+,x^-,h\}$ satisfying
\begin{equation*} [x^+, x^-] =h,\ \ [h,x^\pm] = \pm 2 x^\pm.\end{equation*}
Let $\frak h=\bc h$ be the Cartan subalgebra of $sl_2$, let $\alpha\in\frak h^*$  
the positive root of $sl_2$, given by $\alpha(h)=2$, and set 
$\omega=\alpha/2$. Let $s:\frak h^*\to\frak h^*$ 
be the simple reflection given by $s(\alpha)=-\alpha$.  

The  extended loop algebra of $sl_2$ is the Lie algebra
\begin{equation*}
\eloopg = sl_2\otimes \bc[t,t^{-1}]\oplus \bc 
d,\end{equation*}
with commutator given by
\begin{equation*}
[d, x\otimes t^r] = rx\otimes t^r,  \ \ \ [x\otimes t^r, y\otimes t^s] 
=[x,y]\otimes t^{r+s}\end{equation*}
for $x,y\in sl_2$, $r,s\in\bz$. The loop algebra $\loopg$ is the subalgebra 
$sl_2\otimes \bc[t,t^{-1}]$ of $\eloopg$. 
%The affine Lie algebra $\eloopg$ of $\frak g$ is  one-dimensional central 
%extension of $\eloopg$.
% Let $c$ denote the central element in $\eloopg$. 
Let $\eh =\frak h\oplus\bc d$. Define $\delta\in(\eh)^*$ by
\begin{equation*}
\delta(\frak h) =0,\ \ \delta(d)=1.
\end{equation*}
Extend $\lambda\in\frak h^*$ to an element of $(\eh)^*$ by setting 
$\lambda(d)=0$. Set $P^e =\bz\omega\oplus\bz\delta$, 
and define $P_+^e$ in the obvious way. We regard $s$ as acting on $(\eh)^*$ by 
setting $s(\delta)=\delta$.

%The set of postive roots of $\eloopg$ is then,
%\begin{equation*}\cal R^+ =\cal R^+(<) \cup \cal R^+(0)\cup \cal 
%R^-(>),\end{equation*}
%where,%
%\begin{align*}
%\cal R^+(<) &=\{-\alpha+ k\delta:\alpha\in R^+, k\in\bz, k>0\},\\
%\cal R^+(0) &=\{k\delta:k\in\bz, k>0\},\\
%\cal R^+(>)&=\{\alpha+k\delta:\alpha\in R^+, k\ge 0\}.\end{align*}

For any $x\in sl_2$, $m\in\bz$, we denote by $x_m$ the element $x\otimes 
t^m\in\eloopg$. Set \begin{equation*} e_1^\pm =x^\pm\otimes 1,\ \ e_0^\pm 
=x^\mp\otimes t^{\pm 1}.\end{equation*} Then, the elements $e_i^\pm$, $i=0,1$, 
and $d$ generate $\eloopg$.

For any Lie algebra $\frak a$, the universal enveloping algebra of $\frak a$ 
is 
denoted by $\bu(\frak a)$.   We set 
\begin{equation*}
\bu(\eloopg) =\ebu, \ \ \bu(\loopg)= \bu,\ \ \bu(\frak g)=\ug. 
\end{equation*}
Let $\bu(<)$ (resp. $\bu(>)$) be the subalgebra of $\bu$ generated by  
the $x^-_m$ (resp. $x^+_m$)  for $m\in \bz$.   Set $\ug(<) 
=\bu(<)\cap\ug$ and define $\ug(>)$ similarly. Finally, let $\bu(0)$ be the 
subalgebra of $\bu$ generated by the $h_m$ for all $m\ne 0$. We have 
\begin{align*}
\ug &=\ug(<)\bu(\frak h)\ug(>),\\
\ebu &=\bu(<)\bu(0)\bu(\eh)\bu(>).
\end{align*}

Now let $q$ be an indeterminate, let $\bK=\bc(q)$ be the field of rational
functions in $q$ with complex coefficients, and let $\ba=\bc[q,q^{-1}]$ be
the subring of Laurent polynomials. For $r,m\in\bn$, $m\ge r$, define
\begin{equation*} 
[m]=\frac{q^m -q^{-m}}{q -q^{-1}},\ \ \ \ [m]! =[m][m-1]\ldots [2][1],\ \ \ \ 
\left[\begin{matrix} m\\ r\end{matrix}\right] 
= \frac{[m]!}{[r]![m-r]!}.
\end{equation*}
Then, $\left[\begin{matrix} m\\r\end{matrix}\right]\in\ba$ for all $m\ge r\ge 
0$.

Let $\bu^e_q$ be the  quantized enveloping algebra over $\bK$  associated to 
$\eloopg$. Thus, $\bu^e_q$ is the quotient of the quantum affine algebra 
obtained by setting the central generator equal to $1$. It follows from 
\cite{Dr}, \cite{B}, \cite{J} that $\bu^e_q$ is the algebra 
with generators $\bx_{r}^{{}\pm{}}$ ($r\in\bz$), $K^{{}\pm 1}$, $\bh_{r}$ 
($r\in 
\bz\backslash\{0\}$), $D^{\pm 1}$, 
and the following defining relations:
\begin{align*}
   KK^{-1} = K^{-1}K
  =1, \ \ &\ \ DD^{-1} =D^{-1}D =1,\ \ DK=KD, \\
 K\bh_{r} =\bh_{r}K,\ \ & K\bx_{r}^\pm K^{-1} = q^{{}\pm 2}\bx_{r}^{{}\pm{}},\\ 
 D\bx_{r}^\pm D^{ -1} =q^r\bx^\pm_{r},\ \ & D\bh_rD^{-1}=q^r\bh_r,\\ 
  [\bh_{r},\bh_{s}]=0,\; \; & [\bh_{r} , \bx_{s}^{{}\pm{}}] =
  \pm\frac1r[2r]\bx_{r+s}^{{}\pm{}},\\ 
 \bx_{r+1}^{{}\pm{}}\bx_{s}^{{}\pm{}} -q^{{}\pm 
2}\bx_{s}^{{}\pm{}}\bx_{r+1}^{{}\pm{}} &=q^{{}\pm
    2}\bx_{r}^{{}\pm{}}\bx_{s+1}^{{}\pm{}}
  -\bx_{s+1}^{{}\pm{}}\bx_{r}^{{}\pm{}},\\ [\bx_{r}^+ ,
  \bx_{s}^-]=& \frac{ \psi_{r+s}^+ -\psi_{r+s}^-}{q - q^{-1}},
  \end{align*} 
where the $\psi_{r}^{{}\pm{}}$ are 
determined by equating powers of $u$ in the formal power series 
$$\sum_{r=0}^{\infty}\psi_{\pm r}^{{}\pm{}}u^{{}\pm r} = K^{{}\pm 1} 
{\text{exp}}\left(\pm(q-q^{-1})\sum_{s=1}^{\infty}\bh_{\pm s} u^{{}\pm 
s}\right).$$ 

\vskip 12pt

Define the $q$-divided powers
\begin{equation*}
(\bx_{k}^\pm)^{(r)} 
=\frac{(\bx_{k}^\pm)^r}{[r]!},
\end{equation*}
for all $k\in\bz$, $r\ge 0$.

Define
\begin{equation*}
\blambda^\pm(u) =\sum_{m=0}^\infty {\mbox{\boldmath 
$\Lambda$}}_{\pm m}u^m ={\text{exp}}\left(-\sum_{k=1}^\infty \frac{\bh_{\pm 
k}}{[k]}u^k\right).
\end{equation*}

The subalgebras $\bu_q$, $\bu_q^{\text{fin}}$, $\bu_q(<)$, $\bu(0)$ etc., are 
defined in the 
obvious way. Let $\bu^e_q(0)$ be the subalgebra of $\bu_q^e$ generated by 
$\bu(0)$, $K^{\pm 1}$ and $D^{\pm 1}$. The following result is a simple 
corollary of the PBW theorem for 
$\bu^e_q$, \cite{B}.

\begin{lem}\label{triangle} 
$\bu^e_q=\bu_q(<)\bu^e_q(0)\bu_q(>)$.\hfill\qedsymbol 
\end{lem}

For any invertible element $x\in\bu_q^e$ and any $r\in\bz$, define
\begin{equation*}\left[\begin{matrix}x\\r\end{matrix}\right] = 
\frac{xq^{r} -x^{-1}q^{-r}}{q-q^{-1}}.\end{equation*}

Let $\bu^e_\ba$ be the $\ba$-subalgebra of $\bu^e_q$ generated by the $K^{\pm 
1}$, $(\bx_{k}^\pm)^{(r)}$  ($k\in\bz$,  $r\ge 0$),  
$D^{\pm 1}$ and  $\left[\begin{matrix}D\\ r\end{matrix}\right]$ ($r\in\bz$). 
Then, \cite{L}, \cite{BCP},
\begin{equation*}
\bu^e_q\cong\bu^e_\ba\otimes_\ba\bK.
\end{equation*}

 Define $\bu_\ba(<)$,  $\bu_\ba(0)$ and $\bu_\ba(>)$ in the obvious way. 
Let $\bu^e_\ba(0)$ be the $\ba$-subalgebra of $\bu_\ba$ generated by 
$\bu_\ba(0)$ and the 
elements $K^{\pm 1}$, $D^{\pm 1}$, $\left[\begin{matrix}K\\ 
r\end{matrix}\right]$ and  $\left[\begin{matrix}D\\ r\end{matrix}\right]$ 
($r\in\bz$).
 The following is proved as in Proposition 2.7 in \cite{BCP}.
\begin{prop}
$\bu^e_\ba =\bu_\ba(<)\bu_\ba(0)\bu_\ba^e(\frak h)\bu_\ba(>)$.
\hfill\qedsymbol
\end{prop}

The next lemma is easily checked.
\begin{lem}{\label{auto}}\hfill

\begin{enumerate}
\item[(i)]  There is a unique $\bc$-linear  anti-automorphism $\Psi$ of 
$\bu_q^e$ such that $\Psi(q) =q^{-1}$ and
\begin{align*}\Psi(K) =K,\ \ &\ \Psi(D) =D,\\ \Psi(x_{r}^\pm) =x_{r}^\pm, \ \ 
&\ 
\ \Psi(h_{r}) = -h_{r},\end{align*}
for all $r\in\bz$.
\item[(ii)] There is a unique $\bK$-algebra automorphism $\Phi$ of $\bu^e_q$ 
such that \begin{equation*} \Phi(\bx_{r}^\pm)=\bx_{-r}^\pm,\ \  
\Phi(\Lambda^\pm(u))=\Lambda^\mp(u).\end{equation*}

\item[(iii)] For $0\ne a\in\bK$,  there exists a $\bK$-algebra automorphism 
$\tau_a$ of $\bu_q^e$ such that
\begin{equation*}
\tau_a(\bx_r^\pm)=a^r\bx_r^\pm,\ \ \tau_a(\bh_r)=a^r\bh_r,\ \ \tau_a(K)=K,\ \ 
\tau_a(D)=D,\end{equation*}
for $r\in\bz$. Moreover, 
\begin{equation*}
\tau_a(\blambda_r)=a^r\blambda_r.
\end{equation*}
 \end{enumerate}
\end{lem}

\medspace

{\section {The modules $W_q(m)$} }
 
In this section, we recall the definition and elementary properties of the 
modules $W_q(\lambda)$ from \cite{CP5}, and state the main theorem of this 
paper.
\medskip
\begin{defn} A $\bu_q^e$-module $V_q$  is said to be of {\it type} 1 if 
\begin{equation*}
V_q=\bigoplus_{\lambda\in P^e}(V_q)_\lambda,
\end{equation*}
where the weight space
\begin{equation*} 
(V_q)_\lambda =\{v\in V_q: K.v =q^{\lambda(h)}v,\ \  D.v=q^{\lambda(d)}v\}.
\end{equation*}
A $\bu_q^e$-module of type 1  is said to be {\it integrable} if the elements 
$\bx_{k}^\pm$ act locally nilpotently on $V_q$ for all $k\in\bz$. 
The analogous definitions for $\bu^e$, $\bu^{\text{fin}}$ and 
$\bu_q^{\text{fin}}$ are clear.\end{defn}

\medskip We shall only be interested in modules of type 1 in this paper. It is 
well known \cite{L} that, if $m\ge 0$, there is a unique irreducible 
$\bu_q^{\text{fin}}$-module $V_q^{\text{fin}}(m)$, of dimension $m+1$, 
generated by a vector $v$ such that
\begin{equation*} 
K.v =q^m v,\ \ x_{0}^+.v =0,\ \ (x_{0}^-)^{m+1}.v =0.
\end{equation*} 
Recall \cite{Lbook} that, if $V_q$ is any integrable $sl_2$-module, then
\begin{equation*} 
{\text{dim}}_{\bK}(V_q)_n =\text{dim}_{\bK}(V_q)_{-n},
\end{equation*}
for all $n\in\bz$. Let $V(m)$ denote the $(m+1)$-dimensional irreducible 
representation of $sl_2$.

Define the following generating series in an indeterminate $u$ with 
coefficients in $\bu_q$:
\begin{align*} 
\tilde{\bX}^-(u)=\sum_{m=-\infty}^\infty \bx_{m}^-u^{m+1},\ \ \  
&\ \ \ \bX^-(u)=\sum_{m=1}^\infty \bx_{m}^-u^m,\\
\bX^+(u)=\sum_{m=0}^\infty \bx_{m}^+u^m,\ \ \ &\ \ \ 
\bX_{0}^-(u) =\sum_{m=0}^\infty \bx_{m}^-u^{m+1},\\
\tilde{\bH}(u)=\sum_{m=-\infty}^\infty \bh_{m}u^{m+1},\ \ \ & \ \ \ 
\blambda^\pm(u)=\sum_{m=0}^\infty \blambda_{\pm m}u^m 
={\text{exp}}\left(-\sum_{k=1}^\infty\frac{\bh_{\pm k}}{[k]}u^k\right).
\end{align*}
Given a power series $f$ in $u$, we let $f_s$ denote the coefficient of $u^s$ 
in $f$.

\medskip

For any integer $m\ge 0$, let $I^e_q(m)$ 
be the left ideal in $\bu^e_q$ generated by the elements 
\begin{align*} \bx_k^+\ \ (k\in\bz),\ \  & K-q^m,\ \  D-1,\\ 
\blambda_{r}\ (|r|>m), 
\ \ & \blambda_{m}\blambda_{-r}-\blambda_{m-r}\  (1\le 
r\le m),\\ 
\left(\tilde \bX_i^-(u)\blambda^+(u)\right)_r
\bu(0)\ \  (r\in\bz),\ \ &  \left(\bX_{0}^-(u)^r\blambda^+(u)\right)_s\bu(0) 
\ \ (r\ge 1,\ |s| > m).
\end{align*} 
The ideal $I_q(m)$ in $\bu_q$ is defined in the obvious way (by omitting $D$ 
from the definition).

Set 
\begin{equation*} 
W_q(m) =\bu^e_q/I_q^e(m) \cong\bu_q/I_q(m). 
\end{equation*}
Clearly, $W_q(m)$ is a left $\bu_q^e$-module and a  right $\bu_q(0)$-module. 
Further, the left and right actions of $\bu_q(0)$ on $W_q(m)$ commute. Let 
$w_m$ denote the image of $1$ in $W_q(m)$. If $I_q(m,0)$ (resp. $I_\ba(m,0)$) 
is the left ideal in $\bu_q(0)$ (resp. $\bu_\ba(0)$) generated by the elements 
$\blambda_{m}$ $(|m|>\lambda(h))$ and 
$\blambda_{\lambda(h)}\blambda_{-m}-\blambda_{\lambda(h)-m}$ $(1\le 
m\le\lambda(h))$, then 
\begin{equation*} 
\bu_q(0).w_m \cong\bu_q(0)/I_q(m,0)\ \ ({\text{resp.}}\ \bu_\ba(0).w_m \cong  
\bu_\ba(0)/I_\ba(m,0))
\end{equation*} 
as $\bu_q(0)$-modules (resp. as $\bu_\ba(0)$-modules). The $\bu^e$-modules 
$W(m)$ are defined in the 
analogous way.

Let $\bu_q(+)$ be the  subalgebra of $\bu_q$ generated by the $\bx_k\pm$   for 
$k\ge 0$. The subalgebras $\bu_\ba(+)$ and $\bu(+)$ of $\bu_\ba$ and $\bu$, 
respectively, are defined in the obvious way. The following proposition was 
proved in \cite{CP5}.
\begin{prop}{\label{propwm}} Let $m\ge 1$.\hfill

\begin{enumerate} 
\item[(i)] We have
\begin{align*}
\bu_q(0)/I_q(m,0)&\cong\bK[\blambda_1,\blambda_2,\cdots 
,\blambda_m,\blambda_m^{-1}],\\
\bu_\ba(m,0)/I_\ba(m,0)&\cong\ba[\blambda_1,\blambda_2,\cdots 
,\blambda_m,\blambda_m^{-1}],
\end{align*}
as algebras over $\bK$ and $\ba$, respectively.
\item[(ii)] The $\bu^e$-module $W_q(m)$ is integrable for all $m\ge 0$.
\item[(iii)] $W_q(m)=\bu_q(+).w_m$. In fact, $W_q(m)$ is spannned over $\bK$ by 
the elements
\begin{equation*} (x_0^-)^{(r_0)}(x_1^-)^{(r_1)}\cdots 
(x_{m-1}^-)^{(r_{m-1})}\bu_q(0).w_m,\end{equation*}
where $r_j\ge 0$, $\sum_jr_j\le m$.\end{enumerate}
Analogous results hold for the $\bu$-modules $W(m)$. \hfill\qedsymbol 
\end{prop}

Let $\cal{P}_m$ be the Laurent polynomial ring in $m$ variables with complex 
coefficients. The symmetric group $\Sigma_m$ acts on it in the obvious way; 
let 
$\cal{P}_m^{\Sigma_m}$ be the ring of symmetric Laurent polynomials.
In view of Proposition \ref{propwm}, we see that 
\begin{equation*}\bu_q(0)/I_q(m,0)\cong \bK\cal{P}_m^{\Sigma_m},
\ \ \bu_\ba/I_\ba(m,0)\cong \ba\cal{P}_m^{\Sigma_m},\end{equation*}
where $\bK\cal{P}_m^{\Sigma_m}$ denotes $\cal{P}_m^{\Sigma_m}\otimes\bK$, etc.

Let $V$ be the two-dimensional irreducible $sl_2$-module with basis $v_0,\ v_1$
such that 
\begin{align*} x^+.v_0=0,\ \ h.v_0=v_0,\ \ x^-.v_0 =v_1,\\
x^+.v_1=v_0,\ \ h.v_1=-v_1,\ \ x^-.v_1 =0.\end{align*}
 Let 
$L(V)=V\otimes\bc[t,t^{-1}]$ be the 
$L(sl_2)$-module defined in the obvious way. Let $T^m(L(V))$ be the $m$-fold 
tensor power of $L(V)$ and let $S^m(L(V))$ be its symmetric part. Then, 
$T^m(L(V))$ is a 
left $\bu$-module and a right $\cal{P}_m$-module, and $S^m(L(V))$ is a left 
$\bu$-module and a right $\cal{P}_m^{\Sigma_m}$-module. 
The following was proved in \cite{CP5}.

\begin{thm} As left $\bu$-modules and right  $\cal{P}_m^{\Sigma_m}$-modules, 
we have
\begin{equation*} 
W(m)\cong S^m(L(V)).
\end{equation*}
In particular, $W(m)$ is a free  $\cal{P}_m^{\Sigma_m}$-module of rank 
$2^m$.\hfill\qedsymbol\end{thm}

Our goal in this paper is to prove an analogue of this result for the 
$W_q(m)$.
To do this, we introduce a suitable quantum analogue of $S^m(L(V))$ by 
using the Hecke algebra and a certain quantum symmetrizer.

\begin{defn}
The {\it Hecke algebra} $\cal{H}_m$ is the associative unital algebra over 
$\bc(q)$ 
generated by elements $T_i$ ($i=1,2,\ldots,m-1$)  with the following defining 
relations:
\begin{align*}
(T_i+1)&(T_i-q^2)=0,\\
T_iT_{i+1}T_i&=T_{i+1}T_iT_{i+1},\\
T_iT_j&=T_jT_i\ \ \ {\text{if}}\ |i-j|>1.\ \end{align*}\end{defn}
\medspace

 Set $L_q(V)=L(V)\otimes \bK$. It is easily checked that the following 
formulas define an action of 
$\bu_q^e$ on $L_q(V)$:
\begin{align}
\label{loop} x_k^\pm.(v_\pm\otimes t^r)=0,\ \ &x_k^\pm.(v_\mp\otimes 
t^r)=v_\pm\otimes 
t^{k+r},\\
\label{loop1} \Psi^+(u).(v_\pm\otimes t^r)&=v_\pm\otimes
\frac{q^{\pm 1}-q^{\mp 1}tu}{1-tu}t^r,\\ 
\label{loop2} \Psi^-(u).(v_\pm\otimes t^r)&=v_\pm\otimes
\frac{q^{\mp 1}-q^{\pm 1}t^{-1}u}{1-t^{-1}u}t^r.
\end{align}

The $m$-fold tensor product $T^m(L_q(V))$ is a left $\bu_q^e$-module (the 
action being given by the comultiplication of $\bu_q$) and a right 
$\cal{P}_m$-module (in the obvious way). Now,  
as a vector space over $\bK$, 
\begin{equation*}
L_q(V)^{\otimes m}\cong V^{\otimes m}\otimes_{\bK}\bK[t_1^{\pm 
1},\ldots,t_m^{\pm 1}],
\end{equation*}
and $\Sigma_m$ acts naturally (on the right) on both 
$V^{\otimes m}$ and $\bK[t_1^{\pm 1},\ldots,t_m^{\pm 1}]$ by permuting 
the variables. If $\bv\in V^{\otimes m}$ and $f\in\bK[t_1^{\pm 
1},\ldots,t_m^{\pm 1}]$, denote the action of $\sigma\in\Sigma_m$ by 
$\bv^\sigma$ and $f^\sigma$, respectively. Let $\sigma_i$ be the transposition 
$(i,i+1)\in\Sigma_m$.

\begin{prop}{\label{commute}}\ (\cite[Section 1.2]{KMS}) The Hecke algebra $\cal 
{H}_m$ acts on $L_q(V)^{\otimes m}$ on the right, the 
action of the generators being given as follows :
\begin{equation*}
(v_{t_1}\otimes\cdots\otimes v_{t_m}\otimes f).T_i=
\begin{cases} 
-q(v_{t_1}\otimes\cdots\otimes v_{t_m})^{\sigma_i}\otimes 
f^{\sigma_i}\\
\ \ \ \ \ -(q^2-1)(v_{t_1}\otimes\cdots\otimes v_{t_m})
\otimes\frac{t_{i+1}f^{\sigma_i}-t_if}{t_i-t_{i+1}}\\  
\ \ \ \ \ \ \ \ \ \ {\text{if}}\ t_i=+,\ t_{i+1}=-,\\
-v_{t_1}\otimes\cdots\otimes v_{t_m}\otimes 
f^{\sigma_i}\\
\ \ \ \ -(q^2-1)(v_{t_1}\otimes\cdots\otimes 
v_{t_m})\otimes\frac{t_{i}(f^{\sigma_i}-f)}{t_i-t_{i+1}}\\ 
\ \ \ \ \ \ \ \ \ \ \ {\text{if}}\ t_i=t_{i+1},\\
-q(v_{t_1}\otimes\cdots\otimes v_{t_m})^{\sigma_i}\otimes 
f^{\sigma_i}\\
\ \ \ \ -(q^2-1)(v_{t_1}\otimes\cdots\otimes 
v_{t_m})\otimes\frac{t_{i}(f^{\sigma_i}-f)}{t_i-t_{i+1}}\\ 
\ \ \ \ \ \ \ \ \ \ \ {\text{if}}\ t_i=-,\ t_{i+1}=+.
\end{cases}
\end{equation*}
Moreover, this action commutes with the left action of $\bu^e_q$ on $L_q(V)$ 
and the right action of 
$\bK\cal{P}_m^{\Sigma_m}$.\hfill\qedsymbol\end{prop}

As is well known, the second and third relations in the definition of ${\cal 
H}_m$ imply that, if $\sigma=\sigma_{i_1}\ldots\sigma_{i_N}$ is a reduced 
expression for $\sigma\in\Sigma_m$, so that $N$ is the length $\ell(\sigma)$, 
the element $T_\sigma=T_{i_1}\ldots T_{i_N}\in{\cal H}_m$ depends only on 
$\sigma$, and is independent of the choice of its reduced expression. We 
define 
the symmetrizing operator
\begin{equation*}
{\cal S}^{(m)}:L_q(V)^{\otimes m}\to L_q(V)^{\otimes m}
\end{equation*}
by
\begin{equation*}
{\cal S}^{(m)}=\frac{1}{[m]!}\sum_{\sigma\in\Sigma_m}
(-q^{-2})^{\ell(\sigma)}T_\sigma.
\end{equation*}

\begin{prop}{\label{symmet}} As left $\bu_q^e$-modules and right 
$\bK\cal{P}_m^{\Sigma_m}$-modules, we have
\begin{equation*}
L_q(V)^{\otimes m}={\text{im}}({\cal S}^{(m)})\oplus{\text{ker}}
({\cal S}^{(m)}).
\end{equation*}
\end{prop}
\begin{pf} It is clear from Proposition \ref{commute} that ${\text{im}}({\cal 
S}^{(m)})$  and  ${\text{ker}}
({\cal S}^{(m)})$ are submodules for both the right and left actions.

 The  following proof is  adapted from that of Proposition 1.1 in \cite{KMS}. 
For each $i=1,\ldots,m-1$, we have a factorization
\begin{equation*}
{\cal S}^{(m)}=
\left(\sum_{\sigma'}(-q^{-2})^{\ell(\sigma')}T_{\sigma'}\right)(1-q^{-2}T_i),
\end{equation*}
where $\sigma'$ ranges over $\Sigma_m/\{1,\sigma_i\}$. From this and the first 
of the defining relations of ${\cal H}_m$, it follows that
\begin{equation*}
{\cal S}^{(m)}(T_i+1)=0.
\end{equation*}
In other words, $T_i$ acts on the right on ${\text{im}}({\cal S}^{(m)})$ as 
multiplication by $-1$. It follows that ${\cal S}^{(m)}$ acts on 
${\text{im}}({\cal S}^{(m)})$ by multiplication by the scalar
\begin{equation*}
\frac{1}{[m]!}\sum_{\sigma\in\Sigma_m}(q^{-2})^{\ell(\sigma)}=
\frac{1}{[m]!}\prod_{l=1}^m\frac{1-q^{-2l}}{1-q^{-2}}=q^{-m(m-1)/2}.
\end{equation*}
Hence,
\begin{equation*}
{\cal S}^{(m)}({\cal S}^{(m)}-q^{-m(m-1)/2})=0,
\end{equation*}
and this implies the proposition.\end{pf}

As in \cite{KMS}, 
define an ordered basis $\{u_m\}_{m\in\bz}$ of $L_q(V)$ by setting
\begin{equation*}
u_{-2r}=v_+\otimes t^r,\ \ \ u_{1-2r}=v_-\otimes t^r.
\end{equation*}
Let $u_{r_1}\otimes_S\cdots\otimes_S u_{r_m}$ be the image of 
$u_{r_1}\otimes\cdots\otimes u_{r_m}$ under the projection of $L_q(V)^{\otimes 
m}$ onto $L_q(V)^{\otimes m}/{\text{ker}}({\cal S}^{(m)})$. By Proposition 
2.3, 
this can be identified with an element, which we also denote by 
$u_{r_1}\otimes_S\cdots\otimes_S u_{r_m}$, in ${\text{im}}({\cal S}^{(m)})$. 

\begin{prop}\label{free}
The set $\{u_{r_1}\otimes_S\cdots\otimes_S u_{r_m}\,:\,r_1\ge\cdots\ge r_m\}$ 
is a basis of the vector space ${\text{im}}({\cal S}^{(m)})$. Further,  
${\text{im}}({\cal S}^{(m)})$ is a free $\bK\cal{P}_m^{\Sigma_m}$-module on 
$2^m$ 
generators.
\end{prop}
\begin{pf} The first statement in proved as in \cite{KMS}, Proposition 1.3. 
As for the second, for any $0\le s\le m$, let  ${\text{im}}({\cal S}^{(m)})_s$
be the subspace spanned by {\hbox{$u_{r_1}\otimes_S\cdots\otimes_S u_{r_m}$,}} 
where exactly $s$ of the $r_i$ are even. This space is naturally isomorphic as 
a 
right $\bK\cal{P}_m^{\Sigma_m}$-module to   
$\bK\cal{P}_m^{\Sigma_s\times \Sigma_{m-s}}$. But this module is 
well-known to be free of rank $\left(\begin{matrix}m\\ s\end{matrix}\right)$. 
\end{pf}

Let $\bw=u_0\otimes_S\cdots\otimes_S u_0$. Then, $\bw$ satisfies 
the defining relations of $W_q(m)$, so we have a map of left $\bu_q^e$-modules 
and right $\bK\cal{P}_m^{\Sigma_m}$-modules 
$\eta_m: W_q(m)\to {\text{im}}({\cal S}^{(m)})$ that takes $w_m$ to $\bw$. The 
main theorem of this paper is
\begin{thm} \label{main} The map $\eta_m$ is an isomorphism. In particular,  
$W_q(m)$ is a free 
$\bK\cal{P}_m^{\Sigma_m}$-module of rank $2^m$.
\end{thm}

The theorem is deduced from the following two lemmas.
\begin{lem}{\label{lemma 2.1}} Let $\frak m$ be any maximal ideal in 
$\bK\cal{P}_m^{\Sigma_m}$, and let $d$ be the degree of the field extension 
$\bK\cal{P}_m^{\Sigma_m}/\frak m$ of $\bK$. Then, 
\begin{equation*}\text{dim}_{\bK} \frac{ W_q(m)}{W_q(m)\frak m} 
=2^md.\end{equation*}
\end{lem}

\begin{lem}{\label{lemma 2.2}} The map $\eta_m$ is surjective.
\end{lem}

We defer the proofs of these lemmas to the next section. Once we have these 
two 
lemmas, the proof of Theorem 2 is completed in exactly the same way as Theorem 
1. We include it here for completeness.

\medskip\noindent{\it Proof of Theorem \ref{main}.}  Let $K$ be the kernel of 
$\eta_m$. Since ${\text{im}}({\cal S}^{(m)})$ is a free, hence 
projective, right $\bK\cal{P}_m^{\Sigma_m}$-module by Proposition  
\ref{free}, it follows that
\begin{equation*} 
W_q(m)= {\text{im}}({\cal S}^{(m)})\oplus K,
\end{equation*}
as right $\bK\cal{P}_m^{\Sigma_m}$-modules.
Let $\frak m$ be any maximal ideal in $\bK\cal{P}_m^{\Sigma_m}$. It 
follows from Lemma \ref{lemma 2.1} and Proposition \ref{free} that
\begin{equation*} 
K/K\frak m =0
\end{equation*}
as vector spaces over $\bK$. Since this holds for all 
maximal ideals $\frak m$, Nakayama's lemma implies that $K=0$, proving the 
theorem.
\hfill\qedsymbol
\medskip

{\section {Proof of lemmas \ref{lemma 2.1} and \ref{lemma 2.2}} }
\medskip\noindent In preparation for the proof of  Lemma \ref{lemma 2.1}, we 
first 
show that the modules in question are finite-dimensional.  Recall that a  
maximal ideal in $\bK{\cal P}_m^{\Sigma_m}$ is defined by an $m$-tuple of points 
$\bpi=(\pi_1,\cdots, \pi_m)$, with $\pi_m\ne 0$, in an algebraic closure 
$\overline\bK$ of $\bK$, i.e., it is the kernel of the homomorphism 
$ev_\bpi: \bK{\cal P}_m^{\Sigma_m}\to \overline\bK$ that sends 
$\blambda_i\to\pi_i$. Let $\bF_\bpi$ be the smallest subfield of $\overline\bK$ 
containing $\bK$ and $\pi_1,\cdots,\pi_m$. Clearly, $\bF_{\bpi}$ is a 
finite-rank $\bu_q(0)$-module.
Set
\begin{equation*} W_q(\bpi) = W_q(m) \otimes_{\bu_q(0)} 
\bF_{\bpi},\end{equation*}
and let $w_\bpi=w_m\otimes 1$. The $\bu$-modules $W(\pi)$ are defined 
similarly (with $\pi\in\bc^m$).

The following lemma is immediate from Proposition \ref{propwm}.
\begin{lem}{\label{bpi}} We have 
\begin{equation*} \bu_q(0).w_\bpi = \bF_\bpi w_\bpi.\end{equation*}  
Further, $W_q(\bpi)$ is  spanned over $\bF_\bpi$ by the elements
\begin{equation*} (x_0^-)^{(r_0)}(x_1^-)^{(r-1)}\cdots 
(x_{m-1}^-)^{(r_{m-1})}\end{equation*}
with $\sum _ir_i\le m$.

In particular, $\text{dim}_\bK W_q(\bpi)< \infty$.\hfill\qedsymbol\end{lem}

 The modules $W_q(m)$ and $W_q(\bpi)$, together with their classical 
analogues, have the following universal properties.

\begin{prop}{\label{univ}} Let $\lambda\in P_e^+$.
\begin{enumerate}
\item[(i)] Let $V_q$ be any integrable $\bu_q^e$-module generated by an 
element 
$v$ of $(V_q)_m$ satisfying $\bu_q(>).v =0$. Then, $V_q$ is a quotient of 
$W_q(m)$.
\item[(ii)] Let $V_q$ be a finite-dimensional quotient $\bu_q$-module 
of $W_q(m)$ and let $v$ be the image of $w_m$ in $V_q$.  Assume that 
${\text{ker}}(ev_{\bpi}). v=0$ for some $\bpi=(\pi_1,\cdots ,\pi_m)$, where 
the  $\pi_i\in\overline\bK$.  Then, $V_q$ is a 
quotient of $W_q(\bpi)$.
\item[(iii)] Let $V_q$ be finite-dimensional $\bu_q$-module generated by an 
element $v\in (V_q)_m$ and such that $\bu_q(>).v=0$ and  ${\text{ker}} 
(ev_{\bpi}).v =0$ for some $\bpi$.  
Then, $V_q$ is a quotient of $W_q(\bpi)$. 

Analogous statements hold in the classical case.
\end{enumerate}
\end{prop}
\begin{pf} This proposition was proved  in \cite{CP5} in the case when 
$\bpi\in \bK^m$. The proof in this case is identical, and follows immediately 
from the defining relations of $W_q(m)$ and $W_q(\bpi)$.
\end{pf}

One can now deduce the following theorem, which classifies the irreducible 
finite-dimensional representations of $\bu_q$ over $\bK$. 
 
\begin{thm}{\label{classify}} Let $\bpi\in\overline{\bK}^m$ be as above. Then, 
$W_q(\bpi)$ has a unique irreducible quotient $\bu_q$-module $V_q(\bpi)$. 
Conversely, any irreducible finite-dimensional $\bu_q$-module is isomorphic to 
$V_q(\bpi)$ for a suitable choice of $\bpi$.
\end{thm}
\begin{pf} To prove that $W_q(\bpi)$ has a unique irreducible quotient, it 
suffices to prove that it has a unique maximal $\bu_q$-submodule. For this, it 
suffices to prove that, if $N$ is any submodule, then
\begin{equation*} N\cap W_q(\bpi)_m = \{0\}.\end{equation*}
Since $W_q(\bpi)_m=\bu_q(0).w_\bpi$ is an irreducible $\bu_q(0)$-module, it 
follows that 
\begin{equation*} N\cap W_q(\bpi)_m \ne \{0\}\implies w_\bpi\in N,
\end{equation*}
and hence that $N=W_q(\bpi)$.
Conversely, if $V$ is any finite-dimensional irreducible module, one can show 
as in \cite{CPqa}, \cite{CP5} that there exists $0\ne v\in V_m$ such that 
$\bu_q(>).v=0$ and that $\blambda_r.v =0$ if $|r|>m$.  This shows that
$V_m$ must  be an irreducible module for $\bk[\blambda_1,\cdots ,\blambda_m, 
\blambda_m^{-1}]$, and the result follows.
\end{pf}

It follows from the preceding discussion that, to prove Lemma  \ref{lemma 
2.1}, we must show that, if $\bF_\bpi$ is an extension of $\bK$ of degree $d$, 
then
\begin{equation}{\label{equiv2.1}} \text{dim}_\bK W_q(\bpi) 
=2^md.\end{equation}

Assume from now  on that we have a fixed finite extension $\bF$ of $\bK$ of 
degree $d$ and an element $\bpi\in\bF^m$ as above.
 Given $0\ne a\in\bK$, and $\bpi\in \bF^m$ where $\bK\subset\bF$, define
\begin{equation*} \bpi_a =(a\bpi_1,a^2\bpi_2,\cdots 
,a^m\bpi_m).\end{equation*}

Given any $\bu_q$-module $M$, and $0\ne a\in\bK$, let $\tau_a^*M$ be the 
$\bu_q$-module obtained by pulling back $M$ through the automorphism $\tau_a$ 
defined in Lemma \ref{auto}. The next lemma is immediate from Proposition 
\ref{univ}.
\begin{lem}{\label {twist}} We have
\begin{equation*}
\tau_a^*W_q(m)\cong W_q(m),\ \ \ \tau_a^*W_q(\bpi)\cong W_q(\bpi_a),
\end{equation*}
where the first isomorphism is one of $\bu_q^e$-modules and the second is an 
isomorphism of $\bu_q$-modules.
\hfill\qedsymbol\end{lem}

Let $\overline\ba$ be the integral closure of $\ba$ in $\bF$.
Fix $a\in \ba$ such that $\bpi_{a}\in\overline \ba^m$.
By Lemma 3.2, to prove \eqref{equiv2.1} it suffices to prove that
\begin{equation*}
{\text{dim}}_\bK W_q(\bpi_{a})=2^md.
\end{equation*}
Let $\bl\supset\bk$ be the smallest subfield of $\bF$  such that 
$\bpi_{a}\in\bl^m$ and let $\tilde\ba$ be the integral closure of $\ba$ in 
$\bl$. 
Then, $\tilde\ba$ is free of rank $d$ as an $\ba$-module and 
\begin{equation*}\bl\cong\tilde\ba\otimes_\ba\bK.\end{equation*} 

In what follows we write $\bpi$ for $\bpi_{a}$. Set
\begin{equation*} W_\ba(\bpi)=\bu_\ba\otimes_{\bu_\ba(0)}\tilde\ba 
w_\bpi.\end{equation*}
By Lemma \ref{bpi}, $W_\ba(\bpi)$ is finitely-generated as an 
$\tilde\ba$-module, and hence as an $\ba$-module. Further, 
\begin{equation*} W_q(\bpi)\cong W_\ba(\bpi)\otimes_\ba \bK\end{equation*}
as vector spaces over $\bK$.
Note, however, that $W_\ba(\bpi)$ is not an $\bu_\ba$-module in general, since 
$\pi_m^{-1}$ need not be in $\tilde\ba$. However, $W_\ba(\bpi)$ is a 
$\bu_\ba(+)$-module and  
\begin{equation*} W_q(\bpi)\cong W_\ba(\bpi)\otimes_\ba \bK,\end{equation*}
as $\bu_q(+)$-modules. 

Set 
\begin{equation*}
\bu_1(+)=\bu_\ba(+)\otimes_\ba\bc_1.
\end{equation*}
This  is essentially the universal enveloping algebra $\bu(+)$ of 
$sl_2\otimes\bc[t]$, and hence 
\begin{equation*}
\overline{W_q(\bpi)}=W_\ba(\bpi)\otimes_\ba\bc_1
\end{equation*}
is a module for $\bu(+)$.

Since
\begin{equation*}
{\text{dim}}_\bK W_q(\bpi)={\text{rank}}_\ba 
W_\ba(\bpi)={\text{dim}}_\bc\overline{W_q(\bpi)},
\end{equation*}
 it suffices to prove that
\begin{equation*}
{\text{dim}}_\bc\overline{W_q(\bpi)}=2^md.
\end{equation*}
Define elements $\Lambda_r\in\bu(+)$ in the same way as the elments 
$\blambda_r$ are defined, replacing $q$ by 1.
\begin{lem}
With the above notation, there exists a filtration
\begin{equation*}
\overline{W_q(\bpi)}=W_1\supset W_2\supset\cdots\supset W_d\supset W_{d+1}= 0
\end{equation*}
such that, for each $i=1,\ldots,d$, $W_i/W_{i+1}$ is generated by a non-zero 
vector $v_i$ such that
\begin{align}
x_r^+.v_i=0,\ \ &(x_r^-)^{m+1}.v_i=0\ \ (r\ge 0),\\
h_0.v_i=mv_i,\ \ &\Lambda_r.v_i=\lambda_{i,r}v_i\ \ (r>0),
\end{align}
where the $\lambda_{i,r}\in\bc$ and $\lambda_{i,r}=0$ for $r>m$.
\end{lem}
\begin{pf} Let $\overline{W_q(\bpi)}_n$ be the eigenspace of $h_0$ acting on 
$\overline{W_q(\bpi)}$ with eigenvalue $n\in\bz$. Of course,
\begin{equation*}
\overline{W_q(\bpi)}=\bigoplus_{n=-m}^m\overline{W_q(\bpi)}_n.
\end{equation*}
We can choose a basis $w_1, w_2,\ldots, w_l$, say, of $\overline{W_q(\bpi)}_m$ 
such that the action of $\Lambda_i$, for $i=1,\ldots,m$, is in upper 
triangular form. Let $W_i$ be the $\bu(+)$-submodule of $\overline{W_q(\frak 
m)}$ 
generated by $\{w_i,w_{i+1},\ldots,w_l\}$. This gives a filtration with the 
stated properties. 
To see that $l=d$, note that $W_\ba(\bpi)_m=\tilde\ba w_m$ is a free 
$\ba$-module of rank $d$, hence
\begin{equation*}
\overline{W_q(\bpi)}_m=W_\ba(\bpi)_m\otimes_\ba\bc_1
\end{equation*}
is a vector space of dimension $d$. 
\hfill\end{pf}

\begin{lem} Let $\pi=1+\sum_{r=1}^n\lambda_ru^r\in\bc[u]$ be a polynomial of 
degree 
$n$, and let $m\ge n$. Let $W_+(\pi,m)$ be the quotient of $\bu(+)$ by the 
left ideal generated by the elements
\begin{equation*}
h-m,\ \ \Lambda_r-\lambda_r,\ \ x_r^+,\ \ (x_r^-)^{m+1},
\end{equation*}
for all $r\ge 0$. Then,
\begin{equation*}
{\text{dim}}_\bc W_+(\pi,m)\le 2^m.
\end{equation*}
\end{lem}
\begin{pf} This is exactly the same as the proof given in [CP5, Sections 3 and 
  6] that 
${\text{dim}}_\bc W(\pi)\le 2^{{\text{deg}}(\pi)}$. We note that the arguments 
used there 
only make use of elements of the subalgebra $\bu(+)$ of $\bu$.
\hfill\end{pf}

It follows immediately from this lemma that
\begin{equation*}
{\text{dim}}_\bc\overline{W_q(\bpi)}\le 2^md.
\end{equation*}
Indeed, each $W_i/W_{i+1}$ in Lemma 3.3 is clearly a quotient of some 
$W_+(\pi,m)$ satisfying the conditions of Lemma 3.4, and so has dimension $\le 
2^m$.

\medskip

We have now proved that
\begin{equation*}
{\text{dim}}_\bK W_q(\bpi)\le 2^md.
\end{equation*}
To prove the reverse inequality, let $\tilde\bF$ be the splitting field of the 
polynomial $1+\sum_{i=1}^m\pi_iu^i$ over $\bF$, say
\begin{equation*}
1+\sum_{i=1}^m\pi_iu^i=\prod_{i=1}^m(1-a_iu),
\end{equation*}
with $a_1,\ldots,a_m\in\tilde\bF$.
Let $V_{\bF}(a_i)$ be a two-dimensional vector space over $\tilde\bF$ with 
basis $\{v_+,v_-\}$, define an action of $\bu_q$ on it by setting $t=a_i$ in 
the formulas in
 (\ref{loop}), (\ref{loop1}) and (\ref{loop2}), and set
\begin{equation*} \tilde W=\bigotimes_{i=1}^m V_{\tilde\bF}(a_i).
\end{equation*}
Clearly, 
\begin{equation*}
{\text{dim}}_\bK\tilde W=2^md\tilde d,
\end{equation*}
where $\tilde d$ is the degree of $\tilde\bF$ over $\bF$. If 
$\{f_1,\ldots,f_{\tilde d}\}$ is a basis of $\tilde\bF$ over $\bF$, and if 
$\tilde w =v_+^{\otimes m}$, then
\begin{equation*}
\tilde W=\bigoplus_{j=1}^{\tilde d}\tilde{W}_j,
\end{equation*}
where $\tilde{W}_j$ is the $\bu_q$-submodule of $\tilde W$ generated by 
$f_j\tilde w$ (see \cite[Proof of 2.5]{nato}). Moreover, the vectors 
$f_j\tilde w$ 
satisfy the defining relations of $W_q(\bpi)$, and so are quotients of 
$W_q(\bpi)$. It follows that 
\begin{equation*}
{\text{dim}}_\bK W_q(\bpi)\ge 2^md.\end{equation*}

The proof of Lemma 2.1 is now complete. \hfill\qedsymbol

\medskip

Turning to Lemma \ref{lemma 2.2}, set
\begin{equation*} L_\ba(V)=V\otimes \ba[t,t^{-1}].\end{equation*}
 Clearly,  $L_\ba(V)$ is  a $\bu_\ba$-module. The map ${\cal 
S}^{(m)}$ takes $L_\ba(V)^{\otimes m}$ into itself;
set 
\begin{equation*}  {\text{im}}({\cal S}^{(m)}) =S_q(m),\ \ 
S_\ba(m)= \cal{S}_q(m)\cap L_\ba(V)^{\otimes m}.\end{equation*}
We have
\begin{equation}{\label{class}}
S_\ba^{(m)}\otimes_\ba\bK\cong S_q^{(m)},\ \ \ 
S_\ba^{(m)}\otimes_\ba\bc_1\cong 
S^m(L(V)).
\end{equation}
The first isomorphism above is clear; the second requires the basis 
constructed in Proposition \ref{free}. 
 The proof of Proposition \ref{symmet} shows that
\begin{equation*}
L_\ba(V)^{\otimes m}={\cal S}_\ba(m) \oplus({\text{ker}}({\cal 
S}_q^{(m)})\cap L_\ba(V)^{\otimes m}).
\end{equation*}

Given $\bpi\in \bF^m$ such that $\pi_i\in\overline\ba$, set
\begin{equation*}{S}_q(\bpi) = {S}_q(m)\otimes_{\bu_q(0)}\bF,\ \ \ S_\ba(\bpi) 
={S}_\ba\otimes_{\bu_\ba(0)} \tilde\ba.\end{equation*}
Then, ${S}_q(\bpi)$ (resp. ${S}_\ba(\bpi)$) is  a $\bu_q$-module (resp. 
$\bu_\ba(+)$-module) and
\begin{equation}
S_q(\bpi)\cong S_\ba(\bpi)\otimes_\ba\bK
\end{equation}
as $\bu_q(+)$-modules. Further, the map $\eta_m: W_q(m)\to S_q(m)$ induces a 
map
$\eta_\bpi: W_q(\bpi)\to S_q(\bpi)$ that takes $W_\ba(\bpi)$ into 
$S_\ba(\bpi)$.  

Set $\overline{\bF}=\bF\otimes_{\ba}\bc_1$.  Let $\overline\bpi: 
\bc[\Lambda_1,\cdots ,\Lambda_m]\to\overline\bF$ be the homomorphism obtained 
by sending $\Lambda_i$ to $\pi_i\otimes 1$ and set
\begin{equation*} S(\overline\bpi) = S^m(L(V))\otimes_{\bu(0)}\overline 
\bF.\end{equation*}
Now, in \cite{CP5} we proved that $S^m(L(V))$ is a free $\bc[\Lambda_1,\cdots 
,\Lambda_m]$-module of rank $2^m$, hence $S(\overline\bpi)$ has dimension 
$2^md$. Further, \cite{CP5},
\begin{equation*} W(\overline\bpi)\cong S(\overline\bpi) 
=\bu(+).v_+^{\otimes m}.\end{equation*}
This shows that the induced map
$\overline{\eta_\bpi} :\overline{W_q(\bpi)}\to \overline{S_q(\bpi)}$ is 
surjective and hence, using Lemma \ref{lemma 2.1}, that it is an isomorphism.

Let $K_q(\bpi)$ be  the kernel of $\eta_\bpi$ and let 
$K_\ba(\bpi)=K_q(\bpi)\cap W_q(\bpi)$. Then, $K_\ba(\bpi)$ is free 
$\ba$-module and 
\begin{equation*} \text{dim}_\bk K_q(\bpi)=\text{rank}_\ba 
K_\ba(\bpi).\end{equation*}
The previous argument shows that
\begin{equation*} \overline{K_q(\bpi)} 
=K_\ba(\bpi)\otimes_\ba\bc_1\end{equation*}
is zero. Hence, $K_q(\bpi)=0$ and the map $\eta_\bpi$ is an isomorphism for 
all $\bpi\in \overline{A}^m$. But now, by twisting with an automorphism 
$\tau_a$ 
for $0\ne a\in\bK$,  we have a commutative diagram
 \begin{eqnarray*}
 W_q(\bpi_a)&\longrightarrow&S_q(\bpi_a)\\
 \downarrow&&\downarrow\\
 W_q(\bpi)&\longrightarrow&S_q(\bpi)
 \end{eqnarray*}
for {\it any} $\bpi\in\bF^m$, in which the vertical maps are isomorphisms of 
$\bu_q(+)$-modules. If $a$ is such that $\bpi_a\in\overline{\ba}^m$, the top 
horizontal map is also an isomorphism, hence so is the bottom horizontal map.
Thus, $W_q(\bpi)\to S_q(\bpi)$ is an isomorphism for {\it all} $\bpi\in\bF^m$. 
It follows from Nakayama's lemma that $\eta_m:W_q(m)\to 
S_q^{(m)}$ is surjective and the proof of Lemma 2.2 is 
complete.\hfill\qedsymbol
 
\medskip

{\section{The general case:  a conjecture}} 
In this section, we indicate to what extent the results of this paper can be 
generalized to the higher rank cases, and then state a conjecture in the 
general case.

Thus, let $\frak g$ be  a finite-dimensional simple Lie algebra of rank $n$ of 
type $A,\ D$ or $E$ and let $\hat\frak g$ be the corresponding untwisted 
affine Lie algebra. Given any dominant integral weight $\lambda$ for $\frak 
g$, one can define an integrable $\bu_q(\hat\frak g)$-module $W_q(\lambda)$ on 
which the centre acts trivially, \cite{CP5}. These modules have a family of 
finite-dimensional quotients $W_q(\bpi)$, where 
$\bpi =(\pi^1,\cdots ,\pi^{n})$ and the $\pi^i\in\overline{\bK}^{\lambda(i)}$. 
The module $W_q(\bpi)$ has a unique irreducible quotient 
$V_q(\bpi)$ and one can prove the analogue of Theorem \ref{classify}. (The 
proofs of these statements are the same as in the $sl_2$ case.)

We make the following

\medskip\noindent
{\bf{Conjecture.}}
For any $\bpi$ as above,
\begin{equation*}\text{dim}_\bK W_q(\bpi)= m_{\lambda},\end{equation*}
where $m_\lambda\in\bn$ is given by
\begin{equation*}
m_\lambda=\prod_{i=1}^n (m_i)^{\lambda_i},\ \ \ m_i = \text{dim}_\bK 
W_q(i),\end{equation*}
and $W_q(i)$ is the finite-dimensional module associated to the 
$n$-tuple $(\pi^1,\cdots, \pi^n)$ with $\pi^j=\{0\}$ if $j\ne i$ and 
$\pi^i=\{1\}$.\hfill\qedsymbol

\medskip 

In the case of $sl_2$, the conjecture is established in this paper.
It follows from the results in \cite{C} that $W_q(i)$ is in fact an 
irreducible $\bu_q(\hat\frak g)$-module and hence \cite{banff} the  values of 
the $m_i$ are actually known. The results of \cite{C} also establish the 
conjecture for all $\bpi$ associated to the fundamental weight $\lambda_i$ of 
$\frak g$, for all $i=1,\cdots ,n$.

Using the results in \cite{VV}, one can show that 
\begin{equation*} \text{dim}_\bK W_q(\bpi)\ge m_\lambda.\end{equation*}
It suffices to prove the reverse inequality in the case when 
the $\pi^i\in\overline\ba^{\lambda(i)}$ for all $i$.  One can prove exactly as 
 in this paper that the $\bu_q(+)$-modules $W_q(\bpi)$ admit an 
$\bu_\ba(+)$-lattice $W_\ba(\bpi)$, so that
\begin{equation*} \text{dim}_\bk W_q(\bpi) = \text{rank}_\ba  W_\ba(\bpi)
=\text{dim}_\bc\overline{W_q(\bpi)}.\end{equation*}
Thus, it suffices to prove the conjecture in the classical case, i.e.,
\begin{equation*} \text{dim}_\bc W(\bpi) = m_\lambda,\end{equation*}
where $m_\lambda$ is defined above.

\end{document}